\documentclass[11pt]{article}
\usepackage{amscd}
\usepackage{amsmath,amsfonts,amssymb,amscd}
\usepackage{indentfirst,graphics,epsfig}
 \usepackage{graphicx}
  \DeclareGraphicsExtensions{.eps,.bmp}
\input{epsf}

\usepackage{enumerate}
\usepackage{calc}
\usepackage{lineno}

\makeatletter \@addtoreset{figure}{section} \makeatother
\makeatletter
\long\def\@makecaption#1#2{%
   \vskip 10\p@
   \setbox\@tempboxa\hbox{{#1}\ \ #2}%
   \ifdim \wd\@tempboxa >\hsize
       {#1}\ \ #2\par
   \else
       \hbox to\hsize{\hfil\box\@tempboxa\hfil}%
   \fi}
\makeatother

\newtheorem{thm}{Theorem}[section]
\newtheorem{cor}[thm]{Corollary}
\newtheorem{lem}[thm]{Lemma}

\newtheorem{con}[thm]{Conjecture}

\newcommand{\qed}{{\hfill\rule{3pt}{7pt}}}
\def\pf{\noindent {\it Proof.} }

\def\qed{\hfill \rule{4pt}{7pt}}

\begin{document}
\title{\bf Typical structure of oriented graphs and digraphs with forbidden blow-up transitive triangle \footnote{Research supported by the Training Program for Outstanding Young Teachers in University of Guangdong Province, China (No.312XCQ14564) and the Natural Science Foundation of Guangdong Province, China (296-GK162004).}}

\author{ Jianxi Liu\\
School of Finance,
Guangdong University of Foreign Studies,
Guangzhou, 510006, PR China\\
Email:liujianxi2001@gmail.com\\
}

\date{}
\maketitle

\begin{abstract}
{\small In this work, we establish an analogue result of the Erd\"os-Stone theorem of weighted digraphs using Regularity Lemma of digraphs. We give a stability result of oriented graphs and digraphs with forbidden blow-up transitive triangle and show that almost all oriented graphs and almost all digraphs with forbidden blow-up transitive triangle are almost bipartite respectively.
}\\[3mm]
\end{abstract}
\section{Introduction}
Given a fixed graph $H$, a graph is called $H$-free if it does not contain a subgraph isomorphic to $H$. Denote by $e(G)$ the size (or number of edges) of graph $G$.
Denote by $ex(n,H)$ the maximum size of $H$-free graphs on $n$ vertices.
In the study history of extremal graph theory, there are two types of important problems: (1) For a given graph $H$, determine or estimate $ex(n,H)$, and describle the (asymptotic) strucure of extremal graphs, as $n\rightarrow \infty$. (2) Determine the typical structure of $H$-free graphs on $n$ vertices, as $n\rightarrow \infty$.
The first problem started in 1941 when Tur\'an determined $ex(K_{r+1},n)\le t_r(n):=e(Tu_r(n))$, where the equality holds only by the {\it Tur\'an graph $Tu_r(n)$} which is formed by partitioning the set of $n$ vertices into $r$-parts of nearly equal size, and connecting two vertices by an edge whenever they belong to two different parts.
In 1946 Erd\"os and Stone \cite{ES} extended the Tur\'an theorem and determined $ex({K_{r+1}^t},n)=t_r(n)+o(n^2)$, where $K_{r+1}^t$ is a $K_{r+1}$ blow-up for some positive integer $t$, i.e., $K_{r+1}^t$ is formed by replacing every vertex $v_i$ of $K_{r+1}$ by an independent set of $t$ vertices and connecting every pair of vertices whenever they belong to different independent sets.

The second problem started in 1976 when Erd\"os, Kleitman and Rothschild \cite{EKR} showed that almost all $K_3$-free graphs are bipartite and asymptotically determined the logarithm of the number of $K_r$-free graphs on $n$ vertices, for every integer $r\ge 3$. This was strengthened by Kolaitis, Pr\"omel and Rothschild \cite{KPR}, who showed that almost all $K_r$-free
graphs are $(r-1)$-partite, for every integer $r \ge 3$. Nowadays there are a vast body of work concerning the maximum number of edges and structure of $H$-free graphs on $n$ vertices (see, e.g. \cite{BBS,BBS2, BBS3, BMS, BMSW, EFR, KPR, OPT, PS2}). And some related results have been proved for hypergraphs recently (see, e.g. \cite{BM, PS}).

However, the corresponding questions for digraphs and oriented graphs are almost all wide open, and are the subject of this paper. We shall give some notions before we start to state some relevant results. Given a digraph $G=(V,E)$, let $f_1(G)$ be the number of pairs $u,v\in V$ such that exactly one of $uv$ and $vu$ is an edge of $G$, and let $f_2(G)$ be the number of pairs $u,v\in V$ such that both $uv$ and $vu$ are edges of $G$ (in this case we call $uv$ as a double edge for convenience). 
For a vertex $v$, let $f_1(v)$ be the number of $u\in V$ such that exactly one of $uv$ and $vu$ is an edge of $G$, and let $f_2(v)$ be the number of $u\in V$ such that $uv$ is a double edge. For $a\in \mathbb{R}$ with $a\ge 1$, the {\it weighted size} of $G$ is defined by $e_a(G):= a\cdot f_2(G)+f_1(G)$. For a vertex $v$, its weight is defined by $e_a(v):=a\cdot f_2(v)+f_1(v).$ This definition allows for a unified approach to extremal problems on oriented graphs and digraphs. Because for a digraph $G$, it contains $4^{f_2(G)}2^{f_1(G)}=2^{e_2(G)}$ labelled sub-digraphs and $3^{f_2(G)}2^{f_1(G)}=2^{e_{\log 3}(G)}$ oriented subgraphs if we set $a=2$ and $a=\log 3$, respectively.

Given a digraph $H$, the {\it weighted Tur\'an number $ex_a(n,H)$} is defined as the maximum weighted size $e_a(G)$ among all $H$-free digraphs $G$ on $n$ vertices. Let $DTu_r(n)$ be the digraph obtained from $Tu_r(n)$ by replacing each edge of $Tu_r(n)$ by a double edge. A {\it tournament} is an orientation of a
complete graph. We denote a transitive tournament on $r$ vertices by $T_r$. Note that $DTu_r(n)$ is $T_{r+1}$-free, so $ex_a(n,T_{r+1})\ge e_a(DTu_r(n))=a\cdot t_r(n)$.

For the first problem, Brown and Harary in \cite{BH} determined the extremal digraphs with maximum edges of order $n$ and not containing the transitive tournament $T_{r+1}$. Recently, K\"{u}hn, Osthus, Townsend and Zhao \cite{KOTZ} extended this result to weighted digraphs.
\begin{lem}\cite{KOTZ}\label{LemMaxWeightTFree}
Let $a\in(\frac{3}{2},2]$ be a real number and let $r,n\in \mathbb{N}$. Then $ex_a(n,T_{r+1})=a\cdot t_r(n)$, and $DTu_r(n)$ is the unique extremal $T_{r+1}$-free digraph on $n$ vertices.
\end{lem}
Note that from Lemma \ref{LemMaxWeightTFree} we can see that any $n$-vertex digraph $G$ with $e_a(G)>a\cdot t_r(n)$ contains $T_{r+1}$ for $a\in (\frac{3}{2},2]$. Together with this observation and the Regularity Lemma of digraphs, we establish an analogue Erd\"os-Stone theorem of weighted digraphs as follows:
\begin{thm}\label{Th1}
For all positive integers $r, t$, every real numbers $a\in(\frac{3}{2},2]$ and $\gamma>0$, there exists an integer $n_0$ such that every digraph $G$ with $n\ge n_0$ vertices and
\[ e_a(G) \ge a\cdot t_{r}(n)+\gamma n^2\] contains $T_{r+1}^t$ as a sub-digraph.
\end{thm}

For the second problem, the only results of the above type for oriented graphs were proved by Balogh, Bollob\'as and Morris \cite{BBM, BBM2} who classified the possible `growth speeds' of oriented graphs with a given property.

In 1998 Cherlin \cite{C} gave a classification of countable homogeneous oriented graphs. He remarked that `the striking work of \cite{KPR} does not appear to go over to the directed case' and made the following conjecture.
\begin{con}(Cherlin)\label{con1}
Almost all $T_3$-free oriented graphs are tripartite. 
\end{con}
K\"{u}hn, Osthus, Townsend and Zhao \cite{KOTZ} verified this conjecture and showed that almost all $T_{r+1}$-free oriented graphs and almost all $T_{r+1}$-free digraphs are $r$-partite.
The second part of this work is to reconfirm and generalize the Conjecture \ref{con1}, we show that almost all $T_{3}^t$-free oriented graphs and almost all $T_{3}^t$-free digraphs are almost bipartite for every integer $ t\ge 1$.
More pricisely, let $f(n,T_{3}^t)$ and $f^*(n,T_{3}^t)$ denote the number of labelled $T_{3}^t$-free oriented graphs and digraphs on $n$ vertices, respectively. We show that
\begin{thm}\label{Th2}
For every $r,t\in \mathbb{N}$ with $r\ge 2,t\ge 1$ and any $\alpha >0$ there exists $\epsilon >0$ such that the following holds for all sufficiently large $n$.\\
(i) All but at most $f(n,T_{3}^t)2^{-\epsilon n^2}$ $T_{3}^t$-free oriented graphs on $n$ vertices can be made bipartite by changing at most $\alpha n^2$ edges.\\
(ii) All but at most $f^*(n,T_{3}^t)2^{-\epsilon n^2}$ $T_{3}^t$-free digraphs on $n$ vertices can be made bipartite by changing at most $\alpha n^2$ edges.
\end{thm}
The rest of the paper is organized as followed. We lay out some notations and set out some useful tools in Section 2. We introduce the Regularity Lemma of digraphs and give the proof of Theorem \ref{Th1} in Section 3. We establish a stability result of digraphs and give a proof of Theorem \ref{Th2} in Section 4 and give some conclusion remarks in Section 5.
\section{Notations and Tools}

A {\it digraph} is a pair $(V,E)$ where $V$ is a set of vertices and $E$
is a set of ordered pairs of distinct vertices in $V$ (note that this
means we do not allow loops or multiple edges in the same direction in
a digraph). An {\it oriented graph} is a digraph with at most one edge
between two vertices, so may be considered as an orientation of a simple undirected graph.
In some proofs, given $a,b\in \mathbb{R}$ with $0<a,b<1$, we will use the notation $a \ll b$ to mean that we can find an increasing function $f$ for which all of the conditions in the proof are satisfied whenever $a\le f(b)$. We assume all graphs, oriented graphs and digraphs to be labelled unless otherwise stated. We also assume all large numbers to be integers, so that may some times omit floors and ceilings for the sake of clarity.

Let $G=(V,E)$ be a digraph, we write $uv$ for the edge
directed from $u$ to $v$. For a vertex $v\in V$, we define the out-neighborhood
of $v$ in $G$ to be $N_G^+:=\{u\in V: vu \in E\}$, and the in-neighborhood of $v$ to be $N_G^-:=\{u\in V: uv \in E\}$. The out-degree
$d^+_G(v)$ and the in-degree $d^-_G(v)$ of $v$ in $G$ are defined by $|N_G^+|$
and $|N_G^-|$, respectively.
We define the neighborhood of $v$ to be $N_G(v):= N_G^- \bigcup N_G^+$ and the intersection of out-neighborhood and in-neighborhood of $v$ to be $N^{\pm}_G(v):= N_G^- \bigcap N_G^+$. We write $\Delta(G),\Delta^+(G)$ and $\Delta^-(G)$ for the maximum of $|N_G(x)|, |N^+_G(x)|$ and $|N^-_G(x)|$ over all vertices $v\in G$, respectively. Define $\Delta^0(G)$ as the maximum of $d^+(v)$ and $d^-(v)$ among all $v\in V$.
Given a vertex set $A$ of $G$, the sub-digraph of $G$ induced by $A$ is denoted by $G[A]$ which is the digraph obtained from $G$ by deleting vertices not in $A$ and all their incident edges. Given two disjoint subsets $A$ and $B$ of vertices of $G$, an $A-B$ edge is an edge $ab$ where $a\in A$ and $b\in B$.
We write $E(A,B)$ for the set of all these edges and put $e_G(A,B):=|E(A,B)|$.
We denote by $(A,B)_G$ the bipartite oriented subgraph of $G$ whose vertex class are $A$ and $B$ and whose edge set is $E(A,B)$. The density of $(A,B)_G$ is defined to be
\[  d_G(A,B):=\frac{e_G(A,B)}{|A||B|}.
\]
Given $\epsilon > 0$, we call $(A,B)_G$ is an $\epsilon$-{\it regular} pair if for all subsets $X\subseteq A$ and $Y\subseteq B$ with $|X|>\epsilon |A|$ and $|Y|>\epsilon |B|$ we have that $|d(X,Y)-d(A,B)|<\epsilon$. Note that $(B,A)$ may not be an $\epsilon$-regular pair since the order matters.

For a positive integer $k$ we write $[k]:=\{1,\ldots,k\}$. For convenience, we drop the subscripts of all notions if they are unambiguous. For undefined terminology and notations we refer the reader to \cite{D}.

We need the following result of forbidden digraphs container of K\"uhn et al. \cite{KOTZ}, which allows us to reduce an asymptotic counting problem to an extremal problem.
Given an oriented graph $H$ with $e(H)\ge 2$, we let
$$ m(H)=\max\limits_{H'\subset H, e(H')>1} \frac{e(H')-1}{v(H')-2}.$$

\begin{thm}(\cite{KOTZ}, Theorem 3.3)\label{HFreeContainer}
 Let $H$ be an oriented graph with $h:=v(H)$ and $e(H)\ge 2$, and let $a\in \mathbb{R}$ with $a\ge 1$. For every $\epsilon>0$, there exists $c>0$ such that for all sufficiently large $N$, there exists a collection $\mathcal{C}$ of digraphs on vertex set $[n]$ with the following properties.\\
(a) For every $H$-free digraph $I$ on $[N]$ there exists $G\in \mathcal{C}$ such that $I\subset G$.\\
(b) Every digraph $G\in \mathcal{C}$ contains at most $\epsilon N^h$ copies of $H$, and $e_a(G)\le ex_a(N,H)+\epsilon N^2$.\\
(c) $\log |\mathcal{C}| \le c N^{2-1/m(H)} \log N$.
\end{thm}
Note that this result is essentially a consequence of a recent and very powerful result of Balogh, Morris and Samotij \cite{BMS} and Saxton and Thomason \cite{ST}, which introduces the notion of hypergraph containers to give an upper bound on the number of independent sets in hypergraphs, and a digraph analogue \cite{KOTZ} of the well-known supersaturation result of Erd\"os and Simonovits \cite{ES}.

\section{The Regularity Lemma and Erd\"os-Stone Theorem of Digraphs}
In this section we give the degree form of the regularity lemma for digraphs. A regularity
lemma for digraphs was proved by Alon and Shapira [3]. The degree form follows
from this in the same way as the undirected version (see [34] for a sketch of the
latter). The interested readers can refer to \cite{KS} for a survey on the Regularity Lemma.

\begin{lem}\cite{AS} (Degree form of the Regularity Lemma of Digraphs). 
For all $\epsilon,M'>0$ there exist $M, n_0$
such that if $G$ is a digraph on $n\ge n_0$ vertices and $d\in [0, 1]$, then there exists a
partition of $V(G)$ into $V_0,\ldots,V_k$ and a spanning subdigraph $G'$ of $G$ satisfying the
following conditions:\\
\indent$(1)$ $M'\leq k \leq M$,\\
\indent$(2)$ $|V_0| \leq \epsilon\cdot n$,\\
\indent$(3)$ $|V_1|=\ldots=|V_{k}|=\ell$,\\
\indent$(4)$ $d^+_{G'}(x)>d^+_{G}(x)-(d+\epsilon)n$ for all vertices $x$ of $G$,\\
\indent$(5)$ $d^-_{G'}(x)>d^-_{G}(x)-(d+\epsilon)n$ for all vertices $x$ of $G$,\\
\indent$(6)$ $G'[V_i]$ is empty for all $i=1,\ldots,k$,\\
\indent$(7)$ the bipartite oriented graph $(V_i,V_j)_{G'}$ is $\epsilon$-regular and has density either $0$ or density at least $d$ for all $1\leq i,j \leq k$ and $i\neq j$.
\end{lem}

We call $V_1,\ldots,V_k$ {\it clusters} and $V_0$ the {\it exceptional set}. The last condition of the lemma says that all pairs of clusters are $\epsilon$-regular in both directions (but possibly with different densities). We call the spanning subdigraph $G'\subseteq G$ in the lemma the {\it pure digraph} with parameters $\epsilon, d, \ell$. Given clusters $V_1,\ldots,V_k$ and a digraph $G'$, the reduced digraph $R$ with parameters $\epsilon, d, \ell$ is the digraph whose vertices are $V_1,\ldots,V_k$ and whose edges are all the $V_i-V_j$ edges in $G'$ that is $\epsilon$-regular and has density at least $d$.

Note that a simple consequence of the $\epsilon$-regular pair $(A,B)$: for any subset $Y\subseteq B$ that is not too small, most vertices of $A$ have about the expected number of out-neighbors in $Y$; and similarly for any subset $X\subseteq A$ that is not too small, most vertices of $B$ have about the expected number of in-neighbors in $X$.
\begin{lem}\label{LemExpectNeighbor}
Let $(A,B)$ be an $\epsilon$-regular pair, of density $d$ say, and $X\subseteq A$ has size $|X|\ge \epsilon|A|$ and $Y\subseteq B$ has size $|Y|\ge \epsilon|B|$. Then all but at most $\epsilon|A|$ of vertices in $A$ each of which has at least $(d-\epsilon)|Y|$ out-neighbors in $Y$ and all but at most $\epsilon|B|$ of vertices in $B$ each of which has at least $(d-\epsilon)|X|$ in-neighbors in $X$.
\end{lem}
\pf Let $A'$ be a vertex set with fewer than $(d-\epsilon)|Y|$ out-neighbors in $Y$. Then $e(A',Y)<|A'|(d-\epsilon)|Y|$, so
\[d(A',Y)=\frac{e(A',Y)}{|A'||Y|}<d-\epsilon=d(A,B)-\epsilon.\]
Since $(A,B)$ is $\epsilon$-regular, this implies that $|A'|<\epsilon|A|$.

Similarly, let $B'$ be a vertex set with fewer than $(d-\epsilon)|X|$ in-neighbors in $X$. Then $e(X,B')<|X|(d-\epsilon)|B'|$, so
\[d(X,B')=\frac{e(X,B')}{|X||B'|}<d-\epsilon=d(X,B)-\epsilon.\]
Since $(A,B)$ is $\epsilon$-regular, this implies that $|B'|<\epsilon|B|$.
\qed

The following lemma says that the blow-up $R^s$ of the reduced digraph $R$ can be found in $G$, provided that $\epsilon$ is small enough and the $V_i$ are large enough.

\begin{lem}\label{LemBlowUp}
For all $d\in(0,1)$ and $\Delta\ge 1$, there exists an $\epsilon_0>0$ such that if $G$ is any digraph, $s$ is an integer and $R$ is a reduced digraph of $G'$, where $G'$ is the pure digraph of $G$ with parameters $\epsilon\le \epsilon_0$, $\ell\ge s/\epsilon_0$ and $d$. For any digraph $H$ with $\Delta(G')\le \Delta$, then
\[H\subseteq R^s  \Rightarrow  H\subseteq G'\subseteq G.\]
\end{lem}

\pf The proof is similar with that of Lemma 7.3.2 in \cite{D}.
Given $d$ and $\Delta$, choose $\epsilon_0<d$ small enough that
\begin{align}\label{eq2.1}
\frac{\Delta+1}{(d-\epsilon_0)^\Delta}\epsilon_0 \le 1;
\end{align}
such a choice is possible, since $\frac{\Delta+1}{(d-\epsilon)^\Delta}\epsilon \rightarrow 0$ as $\epsilon \rightarrow 0$. Now let $G,H,s,R$ be given as stated. Let $\{V_0,V_1,\ldots,V_k\}$ be the $\epsilon$-regular partition of $G'$ that give rise to $R$; thus, $\epsilon <\epsilon_0,V(R)=\{V_1,\ldots,V_k\}$ and $|V_1|=\ldots=|V_k|=\ell$. Let us assume that $H$ is actually a sub-digraph of $R^s$, with vertices $u_1,\ldots,u_h$ say. Each vertex $u_i$ lies in one of the $s$-sets $V_j^s$ of $R^s$; this defines a map $\sigma:i\mapsto j$. We aim to define an embedding $u_i \mapsto v_i\in V_{\sigma(i)}$ of $H$ in $G'$; thus, $v_1,\ldots, v_h$ will be distinct, and $v_iv_j$ will be an edge of $G'$ whenever $u_iu_j$ is an edge of $H$.

We choose the vertices $v_1,\ldots,v_h$ inductively. Throughout the induction, we shall have a ``target set" $Y_i \subseteq V_{\sigma(i)}$ assigned to each $i$; this contains the vertices that are still candidates for the choice of $v_i$. Initially, $Y_i$ is the entire set $V_{\sigma(i)}$. As the embedding proceeds, $Y_i$ will get smaller and smaller (until it collapses to $\{v_i\}$): whenever we choose a vertex $v_j$ with $j<i$ and if\\
Case (i): $u_i$ are both out-neighbor and in-neighbor of $u_j$ in $H$, we delete all those vertices from $Y_i$ that are not adjacent to $v_j$ with double edges.\\
Case (ii): $u_i$ is just out-neighbor of $u_j$ in $H$, we delete all those vertices from $Y_i$ that are not the out-neighbor of $v_j$.\\
Case (iii): $u_i$ is just in-neighbor of $u_j$ in $H$, we delete all those vertices from $Y_i$ that are not the in-neighbor of $v_j$.

In order to make this approach work, we have to ensure that the target set $Y_i$ do not get too small. When we come to embed a vertex $u_j$, we consider all the indices $i>j$ such that $u_i$ is adjacent to $u_j$ in $H$; there are at most $\Delta$ such $i$. For each of these $i$, we wish to select $v_j$ so that
\begin{align}\label{eq2.2}
Y_i^j=N^*(v_j)\bigcap Y_i^{j-1}
\end{align}
is large, where
\begin{align*} N^*(v_j)=
\left\{ \begin{array}{ll}
         N^{\pm}(v_j) & \mbox{if $u_i$ are both out-neighbor and in-neighbor of $u_j$};\\
         N^{+}(v_j) & \mbox{if $u_i$ is out-neighbor of $u_j$};\\
         N^{-}(v_j) & \mbox{if $u_i$ is in-neighbor of $u_j$}.\\
          \end{array} \right.
\end{align*}

Now this can be done by Lemma \ref{LemExpectNeighbor}: unless $Y^{j-1}_i$ is tiny (of size less than $\epsilon \ell$), all but at most $\epsilon \ell$ choices of $v_j$ will be such that \eqref{eq2.2} implies
\begin{align}\label{eq2.3}
|Y^j_i| \ge (d-\epsilon)|Y_i^{j-1}|
\end{align}
Doing this simultaneously for all of at most $\Delta$ values of $i$ considered, we find that all but at most $\Delta \epsilon \ell$ choices of $v_j$ from $V_{\sigma(j)}$, and in particular from $Y_j^{j-1}\subseteq V_{\sigma(j)}$, satisfy \eqref{eq2.3} for all $i$.

It remains to show that $|Y^{j-1}|-\Delta \epsilon \ell \ge s$ to ensure that a suitable choice for $v_j$ exists: since $\sigma(j')=\sigma(j)$ for at most $s-1$ of the vertices $u_{j'}$ with $j'<j$, a choice between $s$ suitable candidates for $v_j$ will suffice to keep $v_j$ distinct from $v_1,\ldots,v_{j-1}$. But all this follows from our choice of $\epsilon_0$. Indeed, the initial target sets $Y^0_i$ have size $\ell$, and each $Y_i$ has vertices deleted from it only when some $v_j$ with $j<i$ and $u_j$ and $u_i$ are adjacent in $H$, which happens at most $\Delta$ times. Thus,
\[|Y_i^j|-\Delta \epsilon \ell\ge (d-\epsilon)^\Delta-\Delta \epsilon \ell
\ge (d-\epsilon_0)^\Delta-\Delta \epsilon_0 \ell \ge \epsilon_0 \ell \ge s\]
whenever $j<i$, so in particular $|Y_i^j|-\Delta \ge \epsilon_0 \ell\ge \epsilon \ell$ and $|Y_j^{j-1}|-\Delta \ge \epsilon \ell\ge s$. \qed

We can now prove Theorem \ref{Th1} using Lemma \ref{LemMaxWeightTFree}, Lemma \ref{LemBlowUp} and the Regularity Lemma of digraphs.\\
{\noindent\bf Proof of Theorem \ref{Th1}.}
Let $d:=\gamma, \Delta=\Delta(K_{r+1}^s)$, then Lemma \ref{LemBlowUp} returns an $\epsilon_0>0$. Assume
\begin{align}\label{ineq3.1}
\epsilon_0 < \gamma/2 <1
\end{align}
Let $M'>1/\gamma$, choose $\epsilon>0$ small enough that $\epsilon \le \epsilon_0$ and $\delta:=(a-1)d-\epsilon-a\epsilon^2/2-a\epsilon >0$.
The Regularity Lemma of digraphs returns an integer $M$.
Assume
\[
n \ge \frac{Ms}{\epsilon_0(1-\epsilon)},
\]
Since $\frac{Ms}{\epsilon_0(1-\epsilon)}\ge M'$. The Regularity Lemma of digraphs provided us with an $\epsilon$-regular partition $\{V_0,V_1,\ldots,V_k\}$ of $G'$, the pure digraph of $G$, with parameters $\epsilon,d,\ell$ and $M'\le k\le M$. That is $|V_1|=\ldots=|V_k|=\ell$ and $|V_0|<\epsilon n$. Then
\begin{align}
n\ge k\ell
\end{align}
\[
\ell=\frac{n-|V_0|}{k}\ge \frac{n-\epsilon n}{M}=n\frac{1-\epsilon}{M}\ge \frac{s}{\epsilon_0}
\]
by the choice of $n$. Let $R$ be the regularity digraph of $G'$ with parameters $\epsilon, \ell, d$ corresponding to the above partition.
Since $\epsilon\le \epsilon_0, \ell \ge s/\epsilon_0$. $R$ satisfies the premise of Lemma \ref{LemBlowUp} and $\Delta(K_{r+1}^s)=\Delta$.
Thus in order to conclude by Lemma \ref{LemBlowUp} that $T^s_{r+1}\subseteq G'$, all that remains to be checked is that $T_{r+1}\subseteq R$.

Our plan was to show $T_{r+1}\subseteq R$ by Lemma \ref{LemMaxWeightTFree}. We thus have to checked that the weight of $R$ is large enough.

First by (4) and (5) of the Regularity Lemma of digraphs, we have
\begin{align}\label{eq4.1}
\parallel G\parallel_a \le \parallel G'\parallel_a +(d+\epsilon)n^2
\end{align}

At most ${|V_0| \choose 2}$ double edges lie inside $V_0$, and at most $|V_0|k\ell\le \epsilon n k\ell$ double edges join $|V_0|$ to other partition sets. The $\epsilon-$regular pairs in $G'$ of $0$ density contribute nothing to the weight of $G'$. Since each edge of $R$ corresponds to at most $\ell^2$ edges of $G'$, we thus have in total
\begin{align*}
\parallel G'\parallel_a \le \frac{1}{2}a\epsilon^2 n^2+a\epsilon n k \ell+\parallel R\parallel_a \ell^2
\end{align*}

This together with \eqref{eq4.1}, for all sufficiently large $n$, we have
\begin{align*}
\parallel R\parallel_a
& \ge k^2 \cdot \frac{a(\frac{r-1}{r}+\gamma)n^2-(d+\epsilon)n^2-\frac{1}{2}a\epsilon^2 n^2-a\epsilon n k \ell}{k^2 \ell^2}\\
&\ge a\frac{r-1}{r}k^2+\delta k^2\\
&= a\cdot t_r(k)+\delta k^2\\
& > a\cdot t_r(k).
\end{align*}
Therefore $T_{r+1}\subseteq R$ by Lemma \ref{LemMaxWeightTFree}, as desired.
\qed

Similar with the Erd\"os-Stone theorem of undirected graphs, the Erd\"os-Stone theorem of digraphs is interesting not only in its own right: it also has a most interesting corollary. For an oriented graph $H$, its chromatic number is defined as the chromatic number of its underlying graph. An oriented graph $H$ with chromatic number $\chi(H)$ is called homogeneous if there is an colouring of its vertices by [$\chi(H)$] such that either $E(V_i,V_j)=\emptyset$ or $E(V_j,V_i)=\emptyset$ for every $1\le i\neq j\le \chi(H)$, where $V_i$ is the vertex set with colour $i$.

Given an acyclic homogeneously oriented graph $H$ and an integer $n$, consider the number $h_n:=ex_a(n,H)/(a{n \choose 2})$: the maximum weighted density that an $n-$vertex digraph can have without containing a copy of $H$.

Theorem \ref{Th1} implies that the limit of $h_n$ as $n\rightarrow\infty$ is determined by a very simple function of a natural invariant of $H$--its chromatic number!

\begin{cor}\label{ES}
For every acyclic homogeneously oriented graph $H$ with at least one edge,
\begin{align*}
\lim_{n\rightarrow \infty} \frac{ex_a(n,H)}{a{n \choose 2}}= \frac{\chi(H)-2}{\chi(H)-1}.
\end{align*}
\end{cor}

Before the proof the Corollary \ref{ES} we need the following lemma.
\begin{lem}\cite{D}
\begin{align*}
\lim_{n\rightarrow \infty} \frac{t_{r-1}(n)}{{n \choose 2}}= \frac{r-2}{r-1}.
\end{align*}
\end{lem}

{\bf\noindent Proof of Corollary \ref{ES}}. Let $r:=\chi(H)$. Since $H$ cannot be coloured with $r-1$ colours, we have $H \nsubseteq DTu_{r-1}(n)$ for all $n\in \mathbb{N}$, and hence
\begin{align*}
a t_{r-1}(n)\le ex_a(n,H).
\end{align*}
On the other hand, $H \subseteq T_r^t$
for all sufficiently large $t$, so
\begin{align*}
ex_a(n,H)\le ex_a(n,T_r^t)
\end{align*}
for all those $t$. Let us fix such an $t$. For every $\epsilon>0$, Theorem \ref{Th1} implies that eventually (i.e. for large enough $n$)
\begin{align*}
ex_a(n,T_r^t)<a t_{r-1}(n)+\epsilon n^2.
\end{align*}
Hence for $n$ large,
\begin{align*}
 \frac{t_{r-1}(n)}{{n \choose 2}}
 &\le \frac{ex_a(n,H)}{a{n \choose 2}}\\
 &\le \frac{ex_a(n,T_r^t)}{a{n \choose 2}}\\
 &<\frac{t_{r-1}(n)}{{n \choose 2}}+\frac{\epsilon n^2}{a{n \choose 2}}\\
 &<\frac{t_{r-1}(n)}{{n \choose 2}}+\frac{2\epsilon}{a(1-1/n)}\\
 &\le \frac{t_{r-1}(n)}{{n \choose 2}}+4\epsilon
\end{align*}
Therefore, since $\frac{t_{r-1}(n)}{{n \choose 2}}$ converges to $\frac{r-2}{r-1}$, so does $\frac{ex_a(n,H)}{a{n \choose 2}}$.
\qed

\section{Stability Theorem of Digraphs and Proof of Theorem \ref{Th2}}
In this section, we establish a stability of digraphs and give a proof of Theorem \ref{Th2}. Firstly, we give the result of stability of $T_{3}^t$-free digraphs.
\begin{thm}(Stability Theorem)\label{Stability0}
 Let $a\in R$
with $3/2 <a \le 2$, and $t$ be positive integer. Then for any $T_{3}^t$-free digraph with
\begin{align*}
e_a(G)= a\bigg{(}\frac{1}{2}+o(1)\bigg{)}\frac{n^2}{2}
\end{align*}
satisfies $G=DTu_{2}(n) \pm o(n^2)$.
\end{thm}

\pf
First of all we can assume that all but $o(n)$ vertices of $G$ have weight at least $\frac{an}{2}\big{(}1+o(1)\big{)}$. For otherwise let $v_1,\ldots, v_k, k=\lfloor \epsilon \cdot n \rfloor$
($\epsilon$ is a small positive number independent of $n$) be the vertices of $G$ each of which has weight less than $\frac{an}{2}\big{(}1-c\big{)}$, where $0<c(\epsilon)<c<1$.
But then we have
\begin{align*}
 e_a(G[v_{k+1},\ldots,v_n])
 &\ge (\frac{a}{2}+o(1))\frac{n^2}{2} - \frac{ank}{2}(1-c)\\
& =\bigg{(}\frac{a}{4} (n^2-2kn+k^2)-\frac{k^2}{4}+\frac{ckn}{2}+o(1)\frac{n^2}{2}\bigg{)}\\
& > \frac{a }{4}(n-k)^2 \big{(}1+\delta(\epsilon,c) \big{)},
\end{align*}
where $\delta(\epsilon,c)>0$.
By Theorem \ref{Th1} we have that $G[v_{k+1},\ldots,v_n]$ and therefore $G$ contains a $T_{3}^t$  which contradicts our assumption.

Let now $v_1,\ldots,v_p, p=\big{(}1+o(1) \big{)}n$ be the vertices of $G$ each of which has weight not less than $\frac{an}{2}\big{(}1+o(1)\big{)}$. Then the weight of each vertex of $G[v_1,\cdots,v_p]$ in ($G[v_1,\cdots,v_p]$)
is at least $ap\big{(}\frac{1}{2}+o(1)\big{)}=an\big{(}\frac{1}{2}+o(1)\big{)}$.
And $e_a(G[v_1,\cdots,v_p])=\frac{ap^2}{2}\big{(}\frac{1}{2}+o(1)\big{)}=
\frac{an^2}{2}\big{(}\frac{1}{2}+o(1)\big{)}$. Thus to prove our theorem it will suffice to show that $G[v_1,\cdots,v_p]=DTu_{2}(p)\pm o(p^2)$.

Thus it is clear that without loss of generality we can assume that every vertex of our $G$ has weight at least $an\big{(}\frac{1}{2}+o(1)\big{)}$. Note that we now no longer have to use the assumption of $e_a(G)=\frac{an^2}{2}\big{(}\frac{1}{2}+o(1)\big{)}$. Since our assumption
that $e_a(v_i) \ge an\big{(}\frac{1}{2}+o(1)\big{)}, i=1,\ldots,n$ and $G$ is $T_{3}^t-$free already implies that $e_a(G)=\frac{an^2}{2}\big{(}\frac{1}{2}+o(1)\big{)}$.

We shall show that if $G$ is $T_3^t$-free digraph with $e_a(G)=\frac{an^2}{2}\big{(}\frac{1}{2}+o(1)\big{)}$ for some fixed $t$, then $G=DTu_2(n) \pm o(n^2)$.

A pair of adjacent vertices $u$ and $v$ is called {\it bad} if it is contained in only $o(n)$ of $T_3$ of $G$, otherwise it is called {\it good}.  We divide the proof according the number of good pairs of vertices.

{\bf Case 1.} If $G$ has at least $\alpha n^2$ good pairs of vertices for some $\alpha>0$.\\ 
Let $e_1,\ldots,e_s, s\ge\alpha n^2$ be the edges each of which are contained in at least $\beta n$ of $T_3$, where $\alpha,\beta >0$. We now deduce from this assumption that
$G$ contains a $T_3^t$. Let $v_1^{(i)},\ldots,v_{r_i}^{(i)}$ be the vertices which form a $T_3$ with $e_i, r_i\ge \beta n, s \ge i\ge 1$. Since there are $2^{r}$ orientations of a star $S_{r+1}$ of ${r+1}$ vertices. 
Therefore there are at least $\beta'n:=\beta n/2^{r}$ vertices of $\{v_{j}^{(i)}, r_i\ge j\ge 1\}$ formed with $e_i$ with homogeneous $T_3$, w.l.o.g., assume $\{v_{j}^{(i)}, r'_i\ge j\ge 1\}, r'_i\ge \beta'n$ connect to both end vertices of $e_i$ in the same way. 
Similarly there are at least $\alpha' n^2:=\alpha n^2/2^{r+1}$ edges of $\{e_i,s\ge i\ge 1\}$ each formed with at least $\beta' n$ vertices with homogeneous $T_3$, the addition divisor of two is because there may be two choices of direction of the edges $\{e_i,s\ge i\ge 1\}$. 
And all those $T_3$ formed with those at least $\alpha' n^2$ edges $e_i'$
are homogeneous.

Form all possible $t$-tuple from those homogeneous vertices $v_{r_i}^{(i)}$. We get at least
\begin{align*}
\sum\limits_{i=1}^{\alpha' n^2} {r_i' \choose t}\ge \sum\limits_{i=1}^{\alpha' n^2} {\beta'n \choose t} \ge \alpha' n^2 \frac{(\beta'n)^t}{3^t t!}
>\alpha' n^2 (\frac{\beta'}{3})^t {n \choose t}
\end{align*}
$t$-tuples. Since the total number of $t$-tuples formed from $n$ elements is ${n \choose t}$, there is a $t$-tuple say $z_1,\ldots,z_t$ which corresponds to at least $\alpha' n^2 (\frac{\beta'}{3})^t$ edges $e_i$. By Theorem \ref{Th1} these edges determine a $T_2^t$ with vertices $x_1,\ldots,x_t;y_1,\ldots,y_t$. Thus finally $G[x_1,\ldots,x_t;y_1,\ldots,y_t;z_1,\ldots,z_t]$ and thus $G$ contains a $T_3^t$ as stated. But by our assumption our $G$ does not contain a $T_3^t$. This contradiction completes this part of proof.

{\bf Case 2.} If $G$ has $o(n^2)$ good pairs of vertices. Let $G'$ obtained from $G$ by deleting all edges between every good pair of vertices. Since $e_a(G)= a\bigg{(}\frac{1}{2}+o(1)\bigg{)}\frac{n^2}{2}$, we have 
\[e_a(G)\ge e_a(G')\ge a\bigg{(}\frac{1}{2}+o(1)\bigg{)}\frac{n^2}{2}-a\cdot o(n^2)= a\bigg{(}\frac{1}{2}+o(1)\bigg{)}\frac{n^2}{2},\]
thus $e_a(G')= a\bigg{(}\frac{1}{2}+o(1)\bigg{)}\frac{n^2}{2}$.
By the same argument as in the beginning of the proof, we may assume that $e_a(v_i) \ge an\big{(}\frac{1}{2}+o(1)\big{)}, (i=1,\ldots,n)$ in $G'$. We divide the proof into two subcases according to whether $G'$ contains double edges or not.\\
{\bf Subcase 2.1.} If $G'$ contains double edge(s). Assume $uv$ is a double edge of $G'$, then $u$ and $v$ connect to $(\frac{1}{2}+o(1))n$ vertices with double edges respectively, such that $N(u)\bigcap N(v)=o(n)$.
For otherwise, $N(u)\bigcap N(v)=\Omega(n)$ since both $u$ and $v$ have weight at least $an\big{(}\frac{1}{2}+o(1)\big{)}$. Therefore $u$ and $v$ would be contained in $\Omega(n)$ of $T_3$'s, contradicting to our assumption that $G'$ contains not any good pairs of vertices. \\
{\bf Claim 1.} Every vertex in $N(u)$ ($N(v)$ resp.) has $o(n)$ neighbors in $N(u)$ ($N(v)$, resp.).
For otherwise, say $w\in N(u)$ has $\Omega(n)$ neighbors in $N(u)$, then $uw$ is contained in $\Omega(n)$ of $T_3$'s and is a good pair of vertices, contradicting our assumption.
   
Thus each vertex $w\in N(u)$($w\in N(v)$ resp.) connects to $n\big{(}\frac{1}{2}+o(1)\big{)}$ in $N(v)$($N(u)$ resp.) with double edges. And $e_a(G[N(u)])=o(n^2)$ and $e_a(G[N(v)])=o(n^2)$, then a simple computation shows that  $G$ differs from $DTu_2(|N(u)|,|N(v)|)$ with the vertex set $\{N(u),N(v)\}$ by $o(n^2)$ edges, and $DTu_2(|N(u)|,|N(v)|)$ differs from $DTu_2(n)$ by $o(n^2)$ edges, which prove our theorem (the remaining $n-|N(u)|-|N(v)|=o(n)$ vertices can be clearly ignored).
   
{\bf Subcase 2.2.} If $G'$ does not contain any double edges. Let $UG'$ be its underlying undirected graph. Since $a\in (\frac{3}{2},2]$, we assume $a=\frac{3}{2}+\epsilon$ for some $\epsilon>0$. Then 
\[e_a(v_i)\ge a\bigg{(}\frac{1}{2}+o(1)\bigg{)}n=(\frac{3}{2}+\epsilon)\bigg{(}\frac{1}{2}+o(1)\bigg{)}n=\frac{3}{4}n+\frac{1}{2}\epsilon n+o(n).\]
Assume $uv\in E(G')$, then
\[|N(u)\bigcap N(v)|\ge 2\big(\frac{3}{4}n+\frac{1}{2}\epsilon n+o(n) \big)-n=\frac{1}{4}n+\epsilon n+o(n).\]
For all vertices but $o(n)$ of $N(u)\bigcap N(v)$,  say $w$, we have $wuvw$ is a directed triangle since $uv$ is a bad edge of $G'$. And $w$ only have $o(n)$ neighbors in $N(u)\bigcap N(v)$, for otherwise $uw$ is a good edge. Thus $w$ should have at least $\frac{3}{4}n+\frac{1}{2}\epsilon n+o(n)$ neighbors in the outside of $N(u)\bigcap N(v)$.
But then the number of vertices in $G'$ is at least 
\begin{align*}
&\frac{3}{4}n+\frac{1}{2}\epsilon n+o(n)+|N(u)\bigcap N(v)|\\
&\ge \frac{3}{4}n+\frac{1}{2}\epsilon n+o(n)+\big(\frac{1}{4}n+\epsilon n+o(n)\big)\\
&=n+\frac{3}{2}\epsilon n+o(n)\\
&>n,
\end{align*} 
which is a contradiction and we thus complete the proof.\qed

In order to keep all symbols consistent, we reshape Theorem \ref{Stability0} as follows:\\
{\noindent \bf Theorem of Stability.}\label{Stability}
 Let $a\in R$
with $3/2 <a \le 2$, $t$ be positive integer. Then for any $\beta>0$ there exists $\gamma>0$ such that the following holds for all sufficiently large $n$. If a digraph $G$ on $n$ vertices is $T_{3}^t$-free and
\begin{align*}
e_a(G)= a\bigg{(}\frac{1}{2}-\gamma\bigg{)}\frac{n^2}{2},
\end{align*}
then $G=DTu_2(n) \pm \beta n^2$.

We need the Digraph Removal Lemma of Alon and Shapira \cite{AS}.
\begin{lem}(Removal Lemma).\label{RemovalLemma} For any fixed digraph $H$ on $h$ vertices, and any $\gamma >0$ there exists $\epsilon'>0$ such that the following holds for all sufficiently large $n$. If a digraph $G$ on $n$ vertices contains at most $\epsilon'n^h$ copies of $H$, then $G$ can be made $H$-free by deleting at most $\gamma n^2$ edges.
\end{lem}

We now ready to show that almost all $T_{3}^t$-free oriented graphs and almost all $T_{3}^t$-free digraphs are almost bipartite.

{\noindent\bf Proof of Theorem \ref{Th2}.} We only prove (i) here; the proof of (ii) is almost identical. Let $a:=\log 3$. Choose $n_0 \in \mathbb{N}$ and $\epsilon, \gamma, \beta>0$ such that $1/n_0 \ll \epsilon \ll \gamma \ll \beta \ll \alpha,1/r.$ Let $\epsilon':=2\epsilon$
and $n\ge n_0$. By Theorem \ref{HFreeContainer} (with $T_3^t$ and $\epsilon$ taking the roles of $H,N$ and $\epsilon$ respectively) there is a collection $\mathcal{C}$ of digraphs on vertex set $[n]$ satisfying properties $(a)-(c)$. In particular, every
$T_3^t$-free oriented graph on vertex set $[n]$ is contained in some digraph $G\in \mathcal{C}$. Let $\mathcal{C}_1$ be the family of all those $G \in \mathcal{C}$ for which $e_{\log 3}(G)\ge ex_{\log 3}(n,T_3^t)-\epsilon'n^2$.
Then the number of $T_3^t$-free oriented graphs not contained in some $G\in \mathcal{C}_1$ is at most
\begin{align*}
 |\mathcal{C}| 2^{ex_{\log 3}(n,T_{3}^t)-\epsilon'n^2} \le 2^{-\epsilon n^2} f(n,T_{3}^t),
\end{align*}
because $|\mathcal{C}| \le 2^{n^{2-\epsilon'}}$ and $f(n,T_{3}^t) \ge 2^{ex_{\log 3}(n,T_{3}^t)}$. Thus it suffices to show that every digraph $G\in \mathcal{C}_1$ satisfies $G=DTu_2(n) \pm \alpha n^2$. By (b), each $G\in \mathcal{C}_1$ contains at most $\epsilon'n^{3t}$ copies of $T^t_3$. Thus by Lemma \ref{RemovalLemma} we obtain a $T^t_{3}$-free digraph $G'$ after deleting at most $\gamma n^2$ edges from $G$. Then $e_{\log 3} (G') \ge ex_{\log 3}(n,T_3^t)-(\epsilon'+\gamma)n^2$. We next apply the Theorem of Stability to $G'$ and derive that $G'=DTu_2(n)\pm \beta n^2$. As a result, the original digraph $G$
satisfies $G=DTu_2(n)\pm (\beta+\gamma)n^2$, hence $G=DTu_2(n)\pm \alpha n^2$ as required.\qed

\section{Concluding Remarks}
K\"{u}hn, Osthus, Townsend and Zhao \cite{KOTZ} also gave exactly structures of $T_{r+1}$-free oriented graphs and digraphs, but the exactly structures of $T_{r+1}^t$-free oriented graphs and digraphs
are still out of reach from us. We believe the exact structures are the
same as those of $T_{r+1}$-free oriented graphs and digraphs. Therefore, we ending this paper with the following conjecture:
\begin{con}
Let $r,t\in \mathbb{N}$ with $r\ge 2,t\ge 1$. Then the following hold.\\
(i) Almost all $T^t_{r+1}$-free oriented graph are $r$-partite.\\
(ii) Almost all $T^t_{r+1}$-free digraph are $r$-partite.
\end{con}


\begin{thebibliography}{20}
\bibitem{AS}
N. Alon, A. Shapira, Testing subgraphs in directed graphs, Journal of Computer and System Sciences \textbf{69}(2004), 354-382.

\bibitem{BBM}
J. Balogh, B. Bollob\'as, R. Morris, Hereditary properties of combinatorial structures: posets and oriented
graphs, J. Graph Theory \textbf{56}(2007), 311-332.

\bibitem{BBM2}
J. Balogh, B. Bollob\'as, R. Morris, Hereditary properties of tournaments, Electronic J. Combin. \textbf{14} (2007).



\bibitem{BBS}
J. Balogh, B. Bollob\'as, M. Simonovits, The number of graphs without forbidden subgraphs, J. Combin. Theory
Series B \textbf{91}(2004), 1-24.

\bibitem{BBS2}
J. Balogh, B. Bollob\'as, M. Simonovits, The typical structure of graphs without given excluded subgraphs,
Random Structures and Algorithms \textbf{34}(2009), 305-318.

\bibitem{BBS3}
J. Balogh, B. Bollob\'as, M. Simonovits, The fine structure of octahedron-free graphs, J. Combin. Theory Series B \textbf{101}(2011), 67-84.

\bibitem{BH}
W.G. Brown, F. Harary, Extremal digraphs, in {\it Combinatorial theory and its applicaions}, Coll. Math. Soc. J. Bolyai. \textbf{4}(1970), 135-198.

\bibitem{BM}
J. Balogh, D. Mubayi, Almost all triangle-free triple systems are tripartite, Combinatorica \textbf{32}(2012), 143-169.

\bibitem{BMS}
J. Balogh, R. Morris, W. Samotij, Independent sets in hypergraphs, J. Amer. Math.Soc. \textbf{28}(2014), 669-709.


\bibitem{BMSW}
J. Balogh, R. Morris, W. Samotij, L. Warnke, The typical structure of sparse $K_{r+1}$-free graphs, Trans. Amer. Math. Soc. \textbf{368}(2016), 6439-6485.

\bibitem{C}
G. Cherlin, The classification of countable homogeneous directed graphs and countable homogeneous $n$-tournaments, AMS Memoir \textbf{131}(1998), AMS, Rhode Island.

\bibitem{D}
R. Diestel, Graph Theory, Graduate Texts in Mathematics Vol.173, 4th edition, Springer-Verlag, Heidelberg, 2010.

\bibitem{E}
P. Erd\"os, On extremal problems of graphs and generalized graphs, Israel J. Math. \textbf{2}(1964), 184-190.


\bibitem{EFR}
P. Erd\"os, P. Frankl, V. R\"odl, The asymptotic number of graphs not containing a fixed subgraph and a problem for hypergraphs having no exponent, Graphs and Combin. \textbf{2}(1986), 113-121.

\bibitem{EKR}
P. Erd\"os, D. Kleitman, B. Rothschild, Asymptotic enumeration of $K_n$-free graphs, in {\it Colloquio Internazionale sulle Teorie Combinatorie}(Rome, 1973), Vol. II, 19-27. Atti dei Convegni Lincei 17, Accad. Naz. Lincei, Rome, 1976.

\bibitem{ES}
P. Erd\'os, M. Sinomovits, Supersaturated graphs and hypergraphs, Combinatorica \textbf{3}(1983), 181-192.

\bibitem{KS}
J. Koml\'os, M. Simonovits, Szemer\'edi's Regularity Lemma and its applications in graph theory, Bolyai Society Mathematical Studies 2, Combinatorics, Paul Erd\"os is Eighty (Vol. 2)(D. Mikl\'os, V.T. S\'os and T. Sz\"onyi eds.), Budapest (1996), 295-352.

\bibitem{KOTZ}
D. K\"{u}hn, D. Osthus, T. Townsend and Y. Zhao, On the structure of oriented graphs and digraphs with forbidden tournaments or cycles, to appear in J. Combin. Theory Series B.

\bibitem{KPR}
Ph. Kolaitis, H. Pr\"omel, B. Rothschild, $K_{\ell+1}$-free graphs: asymptotic structure and a $0-1$ law, Trans. Amer. Math. Soc. \textbf{303}(1987), 637-671.

\bibitem{OPT}
D. Osthus, H.J. Pr\"omel, A. Taraz, For which densities are random triangle-free graphs almost surely bipartite? Combinatorica \textbf{23}(2003), 105-150.


\bibitem{PS}
Y. Person, M. Schacht, Almost all hypergraphs without Fano planes are bipartite, in Proceedings of the Twentieth Annual ACM-SIAM Symposium on Discrete Algorithms (SODA, 2009), 217-226. ACM Press.

\bibitem{PS2}
H. Pr\"omel, A. Steger, The asymptotic number of graphs not
containing a fixed color-critical subgraph, Combinatorica \textbf{12}(1992), 463-473.

\bibitem{ST}
D. Saxton, A. Thomason, Hypergraph containers, Invent. Math. \textbf{201}(2015), 925-992.

\end{thebibliography}
\end{document}